\documentclass[12pt]{amsart}
\usepackage{amsmath}
\usepackage{amsthm}
\usepackage{amssymb}
\usepackage{latexsym}

\begin{document}

\newcommand{\curl}{\operatorname{curl}}
\newcommand{\diag}{\operatorname{diag}}
\newcommand{\supp}{\operatorname{supp}}
\newcommand{\comp}{\operatorname{comp}}
\newcommand{\im}{\operatorname{Im}}
\newcommand{\ess}{\operatorname{ess}}
\newcommand{\bu}{{\mathbf u}}
\newcommand{\bp}{{\mathbf p}}
\newcommand{\bq}{{\mathbf q}}
\newcommand{\bk}{{\mathbf k}}
\renewcommand{\Re}{\operatorname{Re}}

\newtheorem{lemma}{Lemma}
\newtheorem{corollary}[lemma]{Corollary}
\newtheorem{proposition}[lemma]{Proposition}
\newtheorem{theorem}[lemma]{Theorem}

\theoremstyle{definition}
\newtheorem{definition}{Definition}[section]

\theoremstyle{remark}

\newtheorem{remark}{Remark}

\title{The Spectrum of a linearized 2D Euler operator}

\author{Y. Latushkin, 
Y. Li, and M. Stanislavova}
\address{Department of Mathematics\\
University of
Missouri\\ Columbia, MO
65211}
\email{yuri@math.missouri.edu}
\email{cli@math.missouri.edu}
\address{Department of Mathematics\\
University of Massachusetts\\ Amherst, MA 01003-4515}
\email{stanis@math.umass.edu}

\begin{abstract}
We study the spectral properties of the linearized  Euler operator
obtained by linearizing  the equations of
incompressible two dimensional fluid at a steady state with the vorticity
that contains only two
nonzero complex conjugate Fourier modes. 
 We prove that the essential spectrum coincides with the imaginary axis,
and give an estimate from above for the number of isolated nonimaginary
eigenvalues. In addition, we prove that the spectral mapping theorem holds for
the group generated by the linearized 2D Euler operator.
\end{abstract}

\maketitle

\section{Introduction}
\label{intro}
In recent years a new interest has been drawn to understanding 
the  {\it stability spectrum} of the Euler equations for 
the motion of inviscid fluid linearized about a steady state. 
We will not even attempt to review a vast literature on the subject, and
refer the readers to the recent survey \cite{Fr} and the bibliography therein,
as well as to the closely related to this paper work in 
\cite{BFY,FH,FVY}.

In this paper, we consider 2D Euler equation
under periodic boundary conditions. We apply Fourier
transform and linearize the Euler equation  
about a steady state that contains
only two nonzero complex conjugate Fourier modes. Two issues are addressed. 

First, we give a 
full description of the spectrum
 of the linearization. Using quite different methods, some
results in this direction were obtained in \cite{L1,L2}.
We show that the 
essential spectrum coincides with the
imaginary axis, and that the discrete nonimaginary spectrum 
consists of finitely many points. Moreover, we give an estimate from 
above for the number of
the  nonimaginary isolated eigenvalues in terms of 
a transparent geometric quantity. Although the existence
of nonimaginary eigenvalues is well-known, see, e.g. \cite{F}, 
we are not aware of any results in the
literature that would give an estimate from above for 
the number of nonimaginary
eigenvalues. The results in \cite{BFY} 
indicate that this estimate is sharp.

Second, we show that the spectral 
mapping theorem holds for the group
generated by the linearized Euler operator $L$, that is, we prove 
that $\sigma(e^{tL})=e^{t\sigma (L)}$, $t\neq 0$, for the spectrum
$\sigma(\cdot)$.
Note that the validity of this spectral mapping
property for (non-analytic) semigroups is a rather delicate issue. 
For instance, the spectral mapping property does not hold even for a
group obtained by a first order perturbation of a two-dimensional wave 
equation, see \cite{Renardy} and many other examples and general 
discussion of this phenomenon in, e.g., \cite{CL,vN}.

Our strategy is to use some ideas from the theory of Jacobi matrices and
operators on spaces with indefinite metrics as well as some general results
from the theory of strongly continuous semigroups. It can be summarized as
follows.

The steady state considered in this paper gives a flow on the torus parallel
to a vector $\bp\in\mathbb{Z}^2$. The linearized Euler operator $L$ is a
difference operator acting on a space of sequences on $\mathbb{Z}^2$, see
\eqref{defLL}. To study the spectrum of $L$ we at first ``slice'' the grid
$\mathbb{Z}^2$ using subsets $\Sigma_\bq$, $\bq\in\mathbb{Z}^2$, of lines
parallel to $\bp$, see \eqref{defSigma}. This gives us a way to represent $L$
as a direct sum of operators $L_\bq$ acting on the space of sequences on
$\mathbb{Z}$. We show that $\sigma(L)=\bigcup_\bq\sigma(L_\bq)$ and
$R(\lambda,L)=\bigoplus_q R(\lambda,L_\bq)$ for the resolvent operators. Using
an appropriate rescaling and ``symmetrization'', we replace the study of
$\sigma(L_\bq)$ by that of $\sigma(B_\bq)$, where $B_\bq$ is a certain 
two-diagonal infinite matrix. The essential spectrum $\sigma_{\ess}(B_\bq)$ is
described using Weyl's Theorem. We show that if $\|\bq\|>\|\bp\|$ then $B_\bq$
is a selfadjoint operator and, therefore, $L_\bq$ does not have nonimaginary
spectrum. If $\|\bq\|\leq \|\bp\|$ then $B_\bq$ is a finite rank perturbation
of a selfadjoint operator. However, we give an appropriate choice of an
indefinite metric that makes $B_\bq$ a $J$-selfadjoint operator on a
Pontryagin space with finitely many positive squares. Standard facts about 
$J$-selfadjoint operators give the estimate for the cardinality of the
nonimaginary point spectrum of $L$ in terms of the number of points in $\mathbb{Z}^2$
located inside of the open disc with the radius $\|\bp\|$. Finally, the proof
of the spectral mapping theorem is based on a general Gearhart-Pr\"{u}ss
theorem, see, e.g., \cite{vN}, and an estimate for the norm of
$R(\lambda,L)$, $\Re\lambda\neq 0$.

\noindent {\bf Acknowledgments.} The authors thank K. Makarov for
many discussions and L. A. Sakhnovich for 
his advice to use $J$--theory. The first author was supported by the Research
Board and Research Council of the University of Missouri. The second author
was supported by the Guggenheim Fellowship and by the Research
Board of the University of Missouri.

\section{Linearization}
\label{sec:1}
Consider the two dimensional Euler equations in the vorticity form:
\begin{equation}\label{Eul} \frac{\partial \Omega }{\partial t}=
-u\frac{\partial \Omega }{\partial x}-
v\frac{\partial \Omega}{\partial y},\quad
\frac{\partial u}{\partial x}+\frac{\partial v}{\partial
y}=0.\end{equation} 
Here $\Omega =\curl \bu=\frac{\partial
v}{\partial x}-\frac{\partial u}{\partial y}$ is the vorticity and
$\bu=(u,v)$ is the velocity. We assume that $u=u(x,y)$ and $v=v(x,y)$ are $2
\pi $-periodic, $x,y\in[0,2\pi]$, and have zero spatial means. 
If $\psi$ is the stream function, then
$u=-\frac{\partial \psi}{\partial y}$, $v=\frac{\partial
\psi}{\partial x}$, and $\Omega =\Delta \psi$. If $\Omega
=\sum_{\bk\in \mathbb{Z}^2\setminus\{0\}}\omega_\bk e^{i\bk\cdot (x,y)}$, 
where $\omega_{-\bk}=\overline{\omega_{\bk}}$, $\bk\in\mathbb{Z}^2\setminus\{0\}$, 
then, after a short calculation, see \cite{L1,L2}, we can rewrite \eqref{Eul} as follows:
\begin{equation}\label{DE} \frac{{d\omega}_\bk}{dt}=
\sum_{\bq\in \mathbb{Z}^2\backslash \{ 0\}}A(\bk-\bq,\bq)
\omega _{\bk-\bq}\omega_\bq,\quad \bk\in
\mathbb{Z}^2.\end{equation}
Here $A(\bp,\bq)$, $\bp,\bq\in \mathbb{Z}^2$, are defined by the formula
$$A(\bp,\bq)=A(\bq,\bp)=\frac{1}{2}
\left[ \frac{1}{\| \bp\|^2}-\frac{1}{\|\bq\|^2}\right]\left|
\begin{array}{cc} p_1 & q_1\\ p_2 & q_2\end{array}\right| ,$$
provided $\bp\neq \pm \bq$, $\bp\neq 0$, $\bq\neq 0$, 
and $A(\bp,\bq)=A(\bq,\bp)=0$
otherwise. Here and below we denote $\bp=(p_1,p_2)$, $\bq=(q_1,q_2)$,
$\bk=(k_1,k_2)$; vertical bars denote the determinant, $\|\cdot\|$ is the
Euclidean norm.

Fix $\bp\in \mathbb{Z}^2\backslash \{0\}$ 
and $\Gamma \in \mathbb{C}$,
and consider the
steady state $\omega^0=(\omega^0_\bk)_{\bk\in\mathbb{Z}^2}$ for \eqref{DE}
defined as follows: 
\begin{equation}\label{defSS}
\omega^0_\bk=\begin{cases}
\frac{1}{2}\Gamma, & \bk=\bp\\ \frac{1}{2}\overline{\Gamma}, & \bk=-\bp
\\ 0, & \bk\neq \pm \bp\end{cases},
\end{equation}
where bar denotes the complex conjugate.

If $\Gamma =a+ib$ and $\bp=(p_1,p_2)$
then the vorticity $\Omega^0$, corresponding to $\omega^0$, is given by the formula
$$\Omega^0(x,y)=a\cos (p_1x+p_2y)-b\sin 
(p_1x+p_2y).$$
The corresponding velocity field has the following components:
\begin{align*} u^0(x,y)&=(-p_2a\sin \bp\cdot (x,y)-p_2b\cos \bp\cdot (x,y))
\| \bp\|^{-2},\\
v^0(x,y)&=(p_1a\sin \bp\cdot (x,y)+p_1b\cos \bp\cdot (x,y))/\|
\bp\|^{-2}.\end{align*} In particular, if $\bp=(0,p_2)$, then the steady
state \eqref{defSS} defines a parallel shear flow 
(see, e.g., \cite{F,Fr}) with
the vorticity $\Omega^0(x,y)=a\cos p_2y-b\sin p_2y$ and the
profile $-ap_2\sin p_2y-bp_2\cos p_2y$. For $b=\im \Gamma=0$,
therefore, we have the sinusoidal profile, studied in \cite{MS,Y}
for the case of a viscous shear flow and, recently, in
\cite{BFY,FSV}. 

The linearization of \eqref{DE} around the steady state
$\omega^0$ as in \eqref{defSS}
gives the
following operator $L$:
\begin{equation}\label{defLL}
L:(\omega_\bk)_{\bk\in \mathbb{Z}^2}\mapsto
\Big(A(\bp,\bk-\bp)\Gamma \omega_{\bk-\bp}+A(-\bp,\bk+\bp)\overline{\Gamma
}\omega_{\bk+\bp}\Big)_{\bk\in \mathbb{Z}^2}.
\end{equation}
The choice of the space on which we want to consider the operator
 $L$ is related to the choice of the space for vorticity in
\eqref{Eul}. Thus, if $\Omega \in H^r(\mathbb{T}^2)$, $r\geq 0$, the 
Sobolev space,
then $L$ should be considered on the space $\ell^2_r(\mathbb{Z}^2)$
of weighted $\ell^2(\mathbb{Z}^2)$-sequences with the weight
$(1+\|\bk\|^2)^{r/2}$, $\bk\in \mathbb{Z}^2$. In what 
follows we will consider the
 operator $L$ on
$\ell^2(\mathbb{Z}^2)$, that is, for $r=0$ and 
$\Omega \in L^2(\mathbb{T}^2)$.

Denote by $W$ the shift operator 
$W:(\omega_\bk)\mapsto (\omega_{\bk-\bp})_{\bk\in \mathbb{Z}^2}$. This is a
unitary operator on $\ell ^2(\mathbb{Z}^2)$.  Also, consider the 
following operator
$D_\bp:(\omega_\bk)\mapsto (A(\bp,\bk)\Gamma \omega_\bk)$, and remark that
$D_{-\bp}=-D_\bp$. Thus, the linearized Euler operator $L$ can be represented as
follows: 
$$L=WD_\bp+W^*D_{-\bp}^*=WD_\bp-W^*D_\bp^*.$$ Here and below
$*$ denotes the adjoint operator, $W^*:(\omega_\bk) \mapsto
(\omega_{\bk+\bp})$ and  $D_\bp^*:(\omega_\bk)\mapsto (A(\bp,\bk) \overline{\Gamma}
\omega_\bk)$.
\begin{remark}\label{rem1}  The definition of $A(\bp,\bq)$ imply that
$D_\bp=D^0_\bp+D^1_\bp$, where
\begin{align*} D^0_\bp&:(\omega_{\bk})\mapsto
 \left(\frac{1}{2\| \bp\|^2}\left|\begin{array}{cc} p_1 & k_1\\p_2
 &k_2\end{array}\right|\Gamma
\omega_{\bk}\right) ,\\
D^1_\bp&:(\omega_{\bk})\mapsto \left( -\frac{1}{2\|
\bk\|^2}\left| \begin{array}{cc} p_1 & k_1\\ p_2 &
k_2\end{array}\right| \Gamma \omega_{\bk}
\right) .\end{align*}
The operator $D^1_\bp$ is a compact operator on
$\ell^2(\mathbb{Z}^2)$.\hfill$\Diamond$\end{remark}

\begin{remark}\label{rem2} Let us define the operator
$L^0=WD^0_\bp-W^*D^{0*}_\bp$. Note that
$$(WD^0_\bp\omega )_{\bk}= {2^{-1}\|
\bp\|^{-2}}\left|\begin{array}{cc} p_1 & k_1-p_1\\ p_2 &
k_2-p_2\end{array}\right| \Gamma \omega_{(
k_1-p_1,k_2-p_2)}$$ implies 
$(WD^0_\bp\omega )_{\bk}=(D^0_\bp W\omega )_{\bk}$.
Therefore, $(L^0)^*=D^{0*}_\bp W^*-D^0_\bp W=-L^0$, and, as a result, $\sigma
(L^0)\subset i\mathbb{R}$. As we will see below,
Weyl's Theorem implies that
$\sigma_{\ess}(L)=\sigma _{\ess}(L^0)$.\hfill$\Diamond$\end{remark}

\section{Decomposition and symmetrization}
\label{sec:2}
In this section, we
 will first perform a ``decomposition'' of the operator $L$. 
The operator $L$ will be
represented as a direct sum of certain operators $L_\bq$ acting
on the space $\ell^2(\mathbb{Z})$. Next, we will perform a ``symmetrization''
of the operators $L_\bq$. This procedure will allow us to study, instead
of the operators $L_\bq$, certain operators $B_\bq$ that are finite rank
perturbations of selfadjoint operators.

{}For $\bp\in \mathbb{Z}^2\backslash \{ 0\}$ as above, and any
$\bq=(q_1,q_2)\in \mathbb{Z}^2$ we denote 
\begin{equation}\label{defSigma}
\Sigma_\bq=\{ \bq+n\bp:n\in
\mathbb{Z}\}.
\end{equation}
 Note that $\Sigma_\bq$ is a subset of the line $\{
\bq+t\bp:t\in \mathbb{R}\}$, but, generally, $\Sigma_\bq$ does not contain
all points with integer coordinates that belong to this line. For
each $\bq\in \mathbb{Z}^2$ select $\widehat{\bq}=\widehat{\bq}(\bq)$ such
that $\widehat{\bq}\in \Sigma_\bq,\| \widehat{\bq}\|=\inf \{ \|
\bq+n\bp\|:n\in \mathbb{Z}\}$, and $\widehat{\bq}=\bq+\max \{ n:\|
\widehat{\bq}\| =\|\bq+n\bp\|\}\bp$. The last condition just takes one of
the two possible points $\bq+n\bp\in \Sigma_\bq$ such that $\|
\bq+n\bp\|=\|\widehat{\bq}\|$. Denote $Q=\{\widehat{\bq}(\bq):\bq\in
\mathbb{Z}^2\}$. Clearly, $\Sigma_\bq=\Sigma_{\widehat{\bq}(\bq)}$ and
$\bigcup_{\widehat{\bq}\in Q}\Sigma_{\widehat{\bq}}=\mathbb{Z}^2$.

{}Fix ${\bq}\in Q$ and let $X_{{\bq}}=\{ (\omega_\bk)\in
\ell^2(\mathbb{Z}^2):\omega_\bk=0$ for $\bk\not\in
\Sigma_{{\bq}}\}$. Note that $\Sigma_{{\bq}}\cap
\Sigma_{{\bq}'}=\emptyset$ provided ${\bq}\neq {\bq}'$;
therefore, $\{ X_{\bq}\}_{{\bq}\in Q}$ is a system of
orthogonal subspaces in $\ell^2(\mathbb{Z}^2)$ such that
$\bigoplus_{{\bq}\in Q}X_{{\bq}}=\ell^2(\mathbb{Z}^2)$.

Note that $J:X_{{\bq}}\to
\ell^2(\mathbb{Z}):(\omega_{{\bq}+n\bp})_{n\in\mathbb{Z}}\mapsto
(\omega_n)_{n\in \mathbb{Z}}$ is an isometric isomorphism. Also,
the operator $L$ leaves the subspace $X_{{\bq}}$ invariant. Thus,
$L=\bigoplus_{{\bq}\in Q}L_{{\bq}}$, where  $L_{{\bq}}=
VD-V^*D^*$, and we denote 
\[V:(\omega_n)\mapsto (\omega_{n-1}),\quad D=D(\bp,{\bq}):(\omega_n)
\mapsto (A({\bq}+n\bp,\bp)\Gamma \omega_n).\]

\begin{remark} If ${\bq}=t\bp$, $t\in \mathbb{R}$, then
$D=0$ by the definition of $A(\bp,\bq)$. Thus, in what follows we assume that
${\bq}\neq t\bp$, that is, that $\Sigma_{{\bq}}$ does not
belong to the line $\{ t\bp:t\in \mathbb{R}\}$. In particular,
${\bq}+n\bp\neq 0$ and ${\bq}+n\bp\neq \pm \bp$ for all $n\in\mathbb{Z}$.
\hfill$\Diamond$\end{remark}

Our next step is to perform a symmetrization of the operator $L_\bq$.
Fix $\bp$ as above, and $\bq\in Q$ such that $\bq\neq t\bp$, 
$t\in \mathbb{R}$. Denote
\begin{equation}\label{defBETA}
\beta =-\frac{1}{2\| \bp\|^2}\left|\begin{array}{cc} q_1 & p_1
\\ q_2 & p_2\end{array}\right|\Gamma,\quad \gamma_n=-\
\frac{\| \bp\|^2}{\|\bq+n\bp\|^2},\quad n\in \mathbb{Z}.\end{equation}
Note that $\beta \neq 0$ due to $\bq\neq t\bp$. By the definition of
$A(\bp,\bq)$ we have $$A(\bq+n\bp,\bp)\Gamma=\frac{1}{2}\left[\frac{1}{\|
\bq+n\bp\|^2}-\frac{1}{\|\bp\|^2}\right] \left|\begin{array}{cc} q_1 &
p_1\\ q_2 & p_2\end{array}\right|\Gamma =\beta (1+\gamma_n).$$
Therefore, \begin{equation}\label{3a}
L_\bq=(V\beta -V^*\overline{\beta})\diag (1+
\gamma_n)_{n\in \mathbb{Z}}.\end{equation} Define $\alpha =i\beta
^{\frac{1}{2}}/\overline{\beta }^{\frac{1}{2}}$, $|\alpha |=1$,
and set $S=\alpha V$. Note that $$V\beta -V^*\overline{\beta
}=-{i}\overline{\beta}^{\frac{1}{2}}
(S+S^*)\beta^{\frac{1}{2}}=-i{|\beta
|}(S+S^*).$$

Denote $M^0_\bq=S+S^*$, observe that $(M^0_\bq)^*=M^0_\bq$, and remark that
\begin{equation}\label{defMZERO} M^0_\bq=\left[\begin{array}{ccccc} \ddots 
& & & &\\ & 0 & \alpha & 0 & \\ 
&\overline{\alpha} & \text{\fbox{$0$}} & \alpha
&\\ &0 & \overline{\alpha} & 0 & \\ & & & & \ddots
\end{array}\right],\quad |\alpha |=1.\end{equation} Let
$M_\bq=M^0_\bq\diag (1+\gamma_n)_{n\in \mathbb{Z}}$, and note that
$L_\bq=-i|\beta|M_\bq$. Thus, \begin{equation}\label{4a}
\sigma (L_\bq)=-i|\beta |\sigma (M_\bq)
\quad\text{and}\quad L=\bigoplus_{\bq\in Q}-i|\beta |M_\bq.\end{equation}
 If $\lambda\not\in \sigma
(L)$, then $\lambda \not\in \sigma (L_\bq)$ and we have $$R(\lambda
,L)=\bigoplus_{\bq\in Q}R(\lambda ,-i|\beta |M_\bq)=\bigoplus_{\bq\in
Q}\frac{i}{|\beta|}R\left( \frac{\lambda}{i|\beta|},M_\bq\right).$$
The operator $M_\bq$ is a multiplicative perturbation of a selfadjoint
operator $M^0_\bq$. Instead of $M_\bq$ we will consider 
a symmetric operator $B_\bq$
with the same spectrum. To do this,
for $n\in \mathbb{Z}$ we define:
\begin{equation}\label{defDELTA} \delta _n=|1+
\gamma_n|^{\frac{1}{2}}\text{ if } 1+\gamma_n\geq 
0\text{ and } \delta_n=i|1+\gamma_n|^{\frac{1}{2}}
\text{
if } 1+\gamma _n\leq 0.\end{equation} Note that
$\delta^2_n=1+\gamma_n$, and $\delta_n\in \mathbb{R}$ if and only
if $\| \bq+n\bp\|\geq \|\bp\|$, and $\delta _n\in i\mathbb{R}$
if and only
if $\| \bq+n\bp\|< \|\bp\|$. Denote
\begin{equation}\label{DefB} B_\bq=\diag (\delta _n)M^0_\bq\diag (\delta _n)
=\left[ \begin{array}{ccccc} \ddots & & & & \\ & 0 & 
\alpha \delta_{-1}\delta_0 &
0 &\\ & \overline{\alpha }\delta_{-1}\delta_0 & \text{\fbox{$0$}} & \alpha 
\delta_0\delta_1 & \\ & 0 & \overline{\alpha} \delta_0\delta_1 & 0 & \ddots
\end{array}\right],\end{equation}
and note that
\begin{equation}\label{MQB} M_\bq=M^0_\bq\diag (\delta_n)\cdot \diag (\delta_
n),\quad B_\bq=\diag (\delta_n)M_\bq^0\diag (\delta_n).\end{equation}

\begin{proposition}\label{PropAB} The
nonzero elements in $\sigma (M_\bq)$ and $\sigma (B_\bq)$ 
coincide.\end{proposition}
\begin{proof}
This is a consequence of the following elementary fact: 
If $A$ and $B$ are bounded operators, 
then $\sigma (AB)\backslash \{0\}=\sigma (BA)\backslash
\{0\}$.\end{proof}

\begin{remark}\label{rem4} \label{qgreaterp} 
Assume $\| \bq\|\geq\|\bp\|$. By the 
choice of $\bq\in Q$ we have $\| \bq+n\bp\|\geq \| \bq\|\geq\|\bp\|$. Therefore, $\delta
_n\in \mathbb{R}$
for all $n\in \mathbb{Z}$ and $B_\bq=B^*_\bq$. 
Hence, $\sigma (B_\bq)\subset \mathbb{R}$ or $\sigma 
(L_\bq)\subset i\mathbb{R}$. Since $L_\bq$ is a bounded
operator for each $\bq\in Q$, we have that $\bigcup_{\| \bq\|\leq \| \bp\|}\sigma
 (L_\bq)$ is a bounded set.
Moreover, $\bigcup_{\bq\in Q}\sigma (L_\bq)\backslash
i\mathbb{R}$ is a bounded set.
\hfill$\Diamond$\end{remark}
\begin{remark} Assume $\|\bq\|<\|\bp\|$. Then $\delta_n\delta_{n+1}\in
i\mathbb{R}$ for exactly {\it two} values of $n\in\mathbb{Z}$. Each of these
values corresponds to a pair of two consecutive points in $\Sigma_\bq$ such
that one of the points lies inside and another outside of the disc with the
radius $\|\bp\|$. For all other values of $n\in\mathbb{Z}$ we have
$\delta_n\delta_{n+1}\in\mathbb{R}$. Thus, $B_\bq$ is a perturbation of a
selfadjoint operator by a rank four skew-selfadjoint operator.
\hfill$\Diamond$\end{remark}

It is quite simple to describe the spectrum of the ``constant coefficients''
infinite matrix $M^0_\bq$.
For $\bq\in Q$ and $\beta=\beta (\bq)$, defined in \eqref{defBETA}, let
\[L^0_\bq=V\beta -V^*\overline{\beta}=-i|\beta 
|M^0_\bq,\] where $M^0_\bq$ is as in
\eqref{defMZERO}. Note that $\| M^0_\bq\|\leq 2$ and, thus, 
$\sigma (M^0_\bq)\subset \{ |z|\leq 2\}$ and $\sigma 
(L^0_\bq)\subset \{|z|\leq |\beta |\}$.

If $\mathbb{F}:\ell^2(\mathbb{Z})\to L^2(\mathbb{T})$ is the
Fourier transform, then $$\mathbb{F}((\lambda
+M^0_\bq)(\omega_n))=(\alpha \overline{z}+\lambda +\overline{\alpha
}z)\mathbb{F} [(\omega_n)](z),\quad |z|=1,$$ for each $\lambda \in
\mathbb{C}$, and $\alpha =i\beta^{\frac{1}{2}}/
\overline{\beta}^{\frac{1}{2}}$. Therefore, $-\lambda \not\in
\sigma (M^0_\bq)$ if and only if the function $(\lambda+M^0_\bq)(\cdot)$
defined by the formula $(\lambda +M^0_\bq)(z)=\alpha
\overline{z}+\lambda +\overline{\alpha }z$ is not equal
to zero for all $|z|=1$. If this is the case, then
$(\lambda+M_\bq^0)^{-1}=\mathbb{F}^{-1}\frac{1}{(\lambda
+M^0_\bq)(z)}\mathbb{F}$. Therefore, the following fact holds.

\begin{proposition}\label{Prop2} \begin{enumerate}\item[(a)] $\sigma (M^0_\bq)=\sigma _{
\ess}(M^0_\bq)=[-2,2]$;
\item[(b)] $\sigma (L^0_\bq)=
\sigma_{\ess}(L^0_\bq)=i[-2|\beta|,2|\beta |]$.
\end{enumerate}
\end{proposition}
We note that very interesting and deep results on the spectrum
of the ``variable coefficients'' infinite matrices of the same type
as $M_\bq$ could be found in
\cite{JN} and the literature on Jacobi matrices cited therein. However, the
results that we need for the specific rate of decay of $\gamma_n$ to zero
for the problem in hand do not seem to be available.

\section{Spectrum}
\label{sec:3}

In this section, we describe the spectrum of the linearized Euler
operator $L$. Let $\sigma_p(\cdot)$ denote the point spectrum. For
$\lambda=a+i\tau$, $a\neq 0$, and $\beta=\beta(\bq)$ 
defined in \eqref{defBETA} we denote $z=\lambda/(i|\beta|)$.
\begin{theorem}\label{LemEsSp}
\begin{enumerate}\item[(a)]
 For each $\bq\in Q$ we have:
$$\sigma_{\ess}(L_\bq)=\sigma (V\beta -V^*{\overline{\beta}})=
i[-2 |\beta |,2 |\beta|].$$
\item[(b)] $\sigma (L)=\bigcup_{\bq\in Q}\sigma (L_\bq)$;
\item[(c)] $\sigma_{\ess}(L)=i\mathbb{R}$ and
$\sigma_p(L)\setminus i\mathbb{R}=\bigcup_{\|\bq\|\le\|\bp\|}
\Big(\sigma_p(L_\bq)\setminus i\mathbb{R}\Big)$ is a bounded set
with accumulation points only on $i\mathbb{R}$.
\end{enumerate}\end{theorem}

\begin{proof} {\bf (a).} Split $L_\bq=L^0_\bq+L^{\comp}_\bq$, where
 $L^0_\bq=V\beta -V^*\overline{\beta}$, and $L^{\comp}_\bq=L^0_\bq\diag 
(\gamma_n)_{n\in \mathbb{Z}}$
is a compact operator due to $\lim_{n\to \infty}|\gamma_n|=0$.
Apply Weyl's theorem (see Lemma XIII.4.3 in \cite{RS4})
for $A=L^0_\bq$ and $B=L_\bq$. Note that $\sigma (A)=i[-2 |\beta |,2
|\beta |]$ has an empty interior in $\mathbb{C}$, and
$\mathbb{C}\backslash \sigma(A)$ consists of only one component.
Since $A-B$ is a compact operator, we have
$\sigma_{\ess}(L^0_\bq)=\sigma _{\ess}(L_\bq)$. 
Using Proposition \ref{Prop2}, we have (a).

{\bf (b)}. Split $L$ into the direct sum of two operators, $L^s$ and 
$L^b$, that correspond to ``small'' and ``big'' values of $\|\bq\|$, that is,
write $L=L^s+L^b$, where $L^s=\bigoplus_{\|\bq\| \leq \|\bp\|}L_
\bq$
and $L^b=\bigoplus_{\|\bq\|>\|\bp\|}L_\bq$. Since $L^s$ is a direct sum of
finitely many operators, we have $$\sigma(L)=(\bigcup_{\| \bq\|\leq
\|\bp\|}\sigma (L_\bq))\bigcup \sigma (L^b),$$ and we need to see only
that $\sigma (L^b)\subset \bigcup_{\|\bq\|>\|\bp\|}\sigma (L_\bq)$. But
$|\beta |=|\beta (\bq)| \to \infty$ as $\| \bq\|\to \infty$. Using (a)
in the theorem, we have: $$\bigcup_{\| \bq\|>\| \bp\|}\sigma (L_\bq)\supset
\bigcup_{\| \bq\|>\|\bp\|}i[-2|\beta |,2 |\beta|]=i\mathbb{R}.$$ Thus, it
suffices to show that $\sigma (L^b)\subset i\mathbb{R}$. Since
$$L^b=\bigoplus_{\| \bq\|>\|\bp\|}L_\bq=\bigoplus_{\|\bq\|>\|\bp\|}-i|\beta
(\bq)|M_\bq,$$ it suffices to prove the following claim:
\begin{equation}\label{FBound} \text{if } \lambda \not \in i\mathbb{R},
\text{ then }
\sup_{\|\bq\|>\|\bp\|}\left\| \frac{1}{\beta (\bq)}
\left( z+M_\bq\right)^{-1}\right\| <\infty.\end{equation} Indeed,
assume that \eqref{FBound} is proved. Then $$\| (\lambda
-L^b)^{-1}\|=\left\|\bigoplus_{\|\bq\|>\|\bp\|}\frac{1}{i|\beta(\bq) |}\left(
z+M_\bq\right)^{-1}\right\|<\infty$$ and
$\lambda \not \in \sigma (L^b)$. 

To prove \eqref{FBound}, fix
$\lambda =a+i\tau $, $a \neq 0$. Note that for any 
$\bq\in Q$ we
have (see \eqref{MQB}):
\begin{equation}\begin{split}\label{9.0} z+M_\bq&=
z+M^0_\bq+M^0_\bq\diag (\gamma _n)\\
&=\left( z+M^0_\bq\right)\left[ I+\left(
z+M^0_\bq\right)^{-1}M^0_\bq\diag
(\gamma_n)\right].\end{split}\end{equation} Note that the resolvent
$(z+M^0_\bq)^{-1}$ in this identity exists since $z\not\in
\mathbb{R}$ and $M^0_\bq$ is a self-adjoint operator; moreover,
\begin{equation}\label{9.1} \left\|
\left( z+M^0_\bq\right)^{-1}\right\|=
\frac{1}{|\text{Im } z|}=\frac{|\beta |}{|a|}.\end{equation}

\begin{proposition}\label{EstBETAGA} There exists a 
constant $c(\bp)$ such that for all $\bq\in Q$, $\bq\neq 0$, we have:
$$|\beta (\bq)|\|\diag (\gamma_n)_{n\in \mathbb{Z}}\| \leq {c(\bp)}/{\| \bq\|
}.$$
\end{proposition}

\begin{proof} Using \eqref{defBETA}, we have:
\begin{equation*}\begin{split} |\beta |\| \diag (\gamma_n)_{n\in \mathbb{Z}
}\| &
=|\beta |\sup_{n\in \mathbb{Z}}|\gamma_n|\leq c\sup_{n\in
\mathbb{Z}}\left|{\left|\begin{array}{cc} q_1 & p_1\\ q_2 &
p_2\end{array}\right|}{\|\bq+n\bp\|^{-2}}\right|\\ 
&=c\sup_{n\in
\mathbb{Z}} \left|{\left|\begin{array}{cc} q_1+np_1 & p_1\\
q_2+np_2 & p_2\end{array}\right|}{\| \bq+n\bp\|^{-2}}\right|\\
& =c\sup_{n\in
\mathbb{Z}}\left|{(\bq+n\bp)\cdot \bp^{\bot}}{\| \bq+n\bp\|^{-2}}\right|,
\end{split}\end{equation*} where $\bp^{\bot}=(p_2,-p_1)$ is the
$\mathbb{Z}^2$-vector, perpendicular to $\bp$. Using Cauchy-Schwartz
inequality, we have that $$|\beta |\| \diag (\gamma _n)_{n\in
\mathbb{Z}}\| \leq c\sup_{n\in \mathbb{Z}}\frac{\| \bq+n\bp\|
\|\bp^{\bot}\|}{\|\bq+n\bp\|^2}\leq \frac{c_1}{\inf_{n\in \mathbb{Z}}\|
\bq+n\bp\|}=\frac{c_1}{\| \bq\|},$$ where the definition of $\bq\in Q$ has
been used. \end{proof}

For $\lambda =a+i\tau$ as above, and $c(\bp)$ from 
Proposition~\ref{EstBETAGA}, fix $\bq_0=\bq_0(a)$ such that 
$\|\bq_0\|>\|\bp\|$ and if $\| \bq\|\geq \|\bq_0\|$
then the inequality
\begin{equation}\label{10.1} \frac{2c(\bp)}{|a|\| \bq\|}\leq \frac{1}{2}
\end{equation}
holds. Note that the set $Q_s:=\{\bq\in Q:\| \bq\|\in [\|\bp\|,\|\bq_0\|]\}$ is
finite, and let $Q_b=\{ \bq\in Q:\|\bq\|>\|\bq_0\|\}$.

Using Proposition \ref{EstBETAGA}, and the inequality $\|
M^0_\bq\|\leq 2$, we have (thanks to \eqref{9.1}) that if $\bq\in
Q_b$, then $$\left\| \left( z+M^0_\bq\right)^{-1}M^0_\bq\diag (\gamma_n)_{n\in
\mathbb{Z}}\right\|\leq \frac{2 |\beta |}{|a|}\| \diag
(\gamma_n)\|\leq \frac{2c(\bp)}{|a|\|\bq\|} \leq \frac{1}{2}.$$ Thus,
the operator $I+(z+M^0_\bq)^{-1}M^0_\bq\diag (\gamma_n)$ is invertible
and, for $\bq\in Q_b$,
\begin{equation}\label{10.2} \left\| \left[ I+\left( z
+M^0_\bq\right)^{-1}M^0_\bq\diag (\gamma_n)\right]^{-1}\right\| \leq
2.\end{equation}

For each $\bq\in Q_s$ we remark that $z\not\in \sigma (B_\bq)$, see
\eqref{MQB}, and, hence, $z\not\in\sigma (M_\bq)$. Since $Q_s$ is
a finite set, we have $$\sup_{\bq\in Q_s}\left\| \frac{1}{|\beta |}\left(
z+M_\bq\right)^{-1}\right\|<\infty.$$

\begin{remark}\label{R10} Let $L^b_s=\bigoplus_{\bq\in Q_s}L_\bq$. Then $L^b_
s$ is a bounded operator. Hence, if $\lambda=a+i\tau$, $a\neq 0$, then
$\| (\lambda-L^b_s)^{-1}\|=O\left( {|\tau |^{-1}}\right)$
as $|\tau |\to \infty$.\hfill$\Diamond$\end{remark}

To finish the proof of \eqref{FBound}, we use \eqref{9.0}-\eqref{9.1} and 
\eqref{10.2}:
\begin{equation}\begin{split}\label{11.1} \sup_{\bq\in Q_b}
\left\| {|\beta |^{-1}}\left( z+M_\bq\right)^{-1}\right\|&\leq \sup_
{\bq\in
Q_b}\frac{2}{|\beta|}\left\| \left( z
+M^0_\bq\right)^{-1}\right\|\\ &\leq \sup_{\bq\in
Q_b}{2}{|\beta|^{-1}}\cdot \frac{|\beta
|}{|a|}=\frac{2}{|a|}.\end{split}\end{equation} This proves
\eqref{FBound} and (b) in the theorem.

{\bf (c).} Since $$\sigma_{\ess}(L)=(\bigcup_{\| \bq\|\leq
\|\bp\|}\sigma_{\ess} (L_\bq))\bigcup \sigma_{\ess} (L^b),$$
the first statement follows from (a) in the theorem. The second
statement follows
from Remark \ref{rem4} and (a) in the theorem.
\end{proof}

Since $L$ is the sum of a skew-adjoint operator $L^0$ and a compact
operator, see Remarks \ref{rem1} and \ref{rem2}, the operator
$L$ generates a strongly continuous group.
The following spectral mapping theorem holds for the group
$\{e^{tL}\}_{t\ge0}$. Its proof is similar to the proof of Theorem 1 in 
\cite{GJLS}.
\begin{theorem}
If $L$ is the linearized Euler operator, then
$$\sigma(e^{tL})=e^{t\sigma(L)}, t\neq 0.$$
\end{theorem}

\begin{proof}
Let $L^b=\bigoplus_{\bq\in Q_b}L_\bq$.
Inequality \eqref{11.1} shows that if $\lambda=a+i\tau$,
$a\neq 0$, then $\|(\lambda-L^b)^{-1}\| \leq 2/|a|$. Using
Remark~\ref{R10} we have that
\begin{equation}\label{estres}
\| (\lambda -L)^{-1}\|=O(1)\quad\text{as}\quad |\tau |\to \infty.
\end{equation}
Now the assertion  follows from the resolvent estimate \eqref{estres}
and the following general 
Gearhart--Pr\"{u}ss spectral mapping theorem: 
On a Hilbert space the spectrum $\sigma(e^{tL})$, $t\neq 0$, 
is the set of the points $e^{\lambda t}$ such that
either $\mu_n=\lambda+2\pi n/t$ belongs to $\sigma(L)$ for some
$n\in\mathbb{Z}$, or the sequence $\{\|R(\mu_n,L)\|\}_{n\in\mathbb{Z}}$
is unbounded, see, e.g. \cite{CL}. Recall, that the spectral mapping property always holds
for the point spectrum. Due to \eqref{estres} we conclude that
$\sigma_{\ess}(e^{tL})$ belongs to the unit circle.
\end{proof}

\section{$J$-Theory}
\label{sec:4}

In this section, we obtain an estimate from above for the 
number of non-imaginary
isolated eigenvalues of the operator $L$.
Recall that $\sigma_p(L_\bq)\backslash i\mathbb{R}$ is empty as soon as $\| \bq
\|\geq \| \bp\|$. Let $\varkappa$ denote the number of points $\bq\in
\mathbb{Z}^2$ that belong to the open disk of radius $\| \bp\|$, 
and such that $\bq\neq t\bp$. Since for such $\bq$ we have 
$(\pm q_1)^2+(\pm
q_2)^2<\| \bp\|^2$, we conclude that $\varkappa$ is even.

\begin{theorem}\label{EstANEV} The number of nonimaginary eigenvalues of $L
$ (counting the multiplicities) does not exceed
$2\varkappa $.\end{theorem}

\begin{proof} Since $L_\bq=0$ for $\bq= t\bp$, only those $L_\bq$ for which $\|
 \bq\|<\|\bp\|$ and $\bq\neq t\bp$ will contribute to the nonimaginary 
point spectrum of
$L$. For each such $\bq$ let $n'_\bq$ denote the 
smallest and $n''_\bq$ the largest integer such that 
$\delta _n\in i\mathbb{R}$ for $n=n'_\bq,\ldots ,n''_\bq$, see
\eqref{defDELTA} for the definition of $\delta_n$. Let $J_\bq=\diag (j_n)_{n\in
 \mathbb{Z}}$, where $j_n=1$ if $n=n'_\bq,\ldots ,n''_\bq$ and $j_n=-1$ 
otherwise. Note
that $J_\bq=J^{-1}_\bq=J^*_\bq$ and
$$J_\bq\diag (\delta_n)=\diag (\delta_n)J_\bq=(\diag (\delta_n))^*.$$
Using \eqref{DefB}, we have that $J_\bq B_\bq J_\bq=B^*_\bq$, that is, 
that the operator $B_\bq$ is $J_\bq$-selfadjoint. Let $\langle \omega ,\omega'\rangle$
be the standard scalar product in $\ell^2(\mathbb{Z})$.
Note that the formula $[\omega ,\omega']=\langle J_\bq\omega ,\omega'\rangle$
defines an indefinite metric on $\ell^2(\mathbb{Z})$. Thus, $\ell^2(\mathbb
{Z})$ is a Pontryagin space with $n''_\bq-n'_\bq$ positive squares. 
By a standard
result of the theory of $J$-selfadjoint operators on Pontryagin spaces
(see, e.g., \cite[Cor. II.3.15]{AI}), we have that the number of nonreal 
eigenvalues of
$B_\bq$ does not exceed $2(n''_\bq-n'_\bq)$. By \eqref{MQB}, \eqref{4a}, 
and 
$\varkappa =\sum_\bq(n''_\bq-n'_\bq)$ we have the result.
\end{proof}

\begin{proposition}\label{SymSp} The nonimaginary eigenvalues 
of $L$ are symmetric about the coordinate axes.\end{proposition}

\begin{proof} Re-write \eqref{3a} as $L_\bq=N\diag (1+\gamma_n)$, where 
$N=V\beta -V^*{\bar{\beta}}=-N^*$. Using the argument in the
proof of Proposition~\ref{PropAB}, we
have:
\begin{equation*}\begin{split} \sigma(L^*_\bq)\backslash 
\{0\}&=\sigma (\diag (1+\gamma_n)N^*)\backslash \{ 0\}=-\sigma (\diag
(1+\gamma_n)N)\backslash \{ 0\}\\
&=-\sigma (N\diag (1+\gamma_n))\backslash \{ 0\}
=-\sigma (L_\bq)\backslash \{0\}.\end{split}\end{equation*}

Thus, $\sigma (L_\bq)\backslash \{ 0\}=\overline{\sigma (L^*_\bq)}\backslash 
\{0\}=-\overline{\sigma(L_\bq)}\backslash \{0\}$, and the nonimaginary
eigenvalues of $L$ are symmetric about $i\mathbb{R}$. Next, define $\hat{J}
$ on $\ell^2(\mathbb{Z})$ by $\hat{J}:(\omega_n)\mapsto ((-1)^n\omega_n)$
and note that $V\hat{J}=-\hat{J}V$ and $\hat{J}V^*=-V^*\hat{J}$. Thus,
$$\hat{J}L_\bq\hat{J}=\hat{J}(V\diag (1+\gamma_n)\beta -V^*\diag (1+\gamma_
n)\bar{\beta})\hat{J}=-L_\bq,$$
and $\sigma (L_\bq)=\sigma (\hat{J} L_\bq\hat{J})=-\sigma (L_\bq)$.
\end{proof}

Note that the change of variables
$\xi=p_1x+p_2y$, $\eta =-p_2x+p_1y$ in \eqref{Eul} converts
$\Omega^0$ to the vorticity
$\widetilde{\Omega}^0=\widetilde{\Omega}^0(\xi)$, that is, to a
parallel shear flow. However, this does not simplify our analysis, since
the new flow looses $2\pi$-periodicity, and the results from, say,
\cite{F} cannot be applied directly.

Recently, a very interesting case of a steady state whose vorticity
has four (symmetric) nonzero Fourier modes was considered in \cite{FVY}.
For general $\bp=(p_1,p_2)$, the steady state
considered in the current paper is in an ``intermediate position''
between the parallel shear flow as in \cite{F} and the Kolmogorov
flow as in \cite{BFY}, and the more sophisticated case studied in
\cite{FVY}. For the case in \cite{FVY}, the continuous spectrum
of the linearization is unstable, while in our case it is
always stable (that is, located on the imaginary axis). Note that
the vorticity $\Omega^0$ for the case considered in \cite{FVY}
has the following representation: $$\Omega^0 (x,y)=\Re
({\Gamma_1}e^{i{\bp_1}\cdot
(x,y)}+{{\Gamma_2}}e^{i{{\bp_2}}\cdot
(x,y)})),\quad \Gamma_{1,2}\in\mathbb{C},\quad \bp_{1,2}\in\mathbb{Z}^2.$$

\end{document}